# A sharp-interface mesoscopic model for dendritic growth

**Mitja Jančič**[1,2], **Miha Založnik**[3], **Gregor Kosec**[1]

[1] Institut Jožef Stefan, department of distributed and parallel systems, Jamova cesta 39, 1000 Ljubljana, Slovenia

[2] International postgraduate school Jožef Stefan, Jamova cesta 39, 1000 Ljubljana, Slovenia

[3] Université de Lorraine, CNRS, IJL, F-54000 Nancy, France

E-mail: Gregor.Kosec@ijs.si

**Abstract**. The grain envelope model (GEM) describes the growth of envelopes of dendritic crystal grains during solidification. Numerically the growing envelopes are usually tracked using an interface capturing method employing a phase field equation on a fixed grid. Such an approach describes the envelope as a diffuse interface, which can lead to numerical artefacts that are possibly detrimental. In this work, we present a sharp-interface formulation of the GEM that eliminates such artefacts and can thus track the envelope with high accuracy. The new formulation uses an adaptive meshless discretization method to solve the diffusion in the liquid around the grains. We use the ability of the meshless method to operate on scattered nodes to accurately describe the interface, i.e., the envelope. The proposed algorithm combines parametric surface reconstruction, meshless discretization of parametric surfaces, global solution construction procedure and partial differential operator approximation using monomials as basis functions. The approach is demonstrated on a two-dimensional $h$-adaptive solution of diffusive growth of dendrites and assessed by comparison to a conventional diffuse-interface fixed-grid GEM.

## 1. Introduction

Dendritic grains are the most common growth morphology in solidification of metallic alloys. Their growth is to a large extent governed by collective interactions with adjacent grains. The Grain Envelope Model (GEM) [1,2] is a method that can simulate large ensembles of grains and thus fully describe the collective interactions. The GEM is based on the idea that we do not need to represent individual branches of the dendrite but rather represent the dendrite by its envelope, i.e., a smooth surface that connects all growing dendrite tips and has a much simpler shape. The growth of the envelope is obtained from an analytical stagnant film model that relates the growth speed to the solute concentration in the vicinity of the envelope. This formulation allows us to use numerical mesh spacings that are around an order of magnitude coarser than for methods that describe the solid-liquid interface (e.g., phase field).

The accuracy of the GEM simulations [2,3] is principally controlled by the resolution of the sharp variations of solute concentration in the vicinity of the evolving envelopes and by the algorithm used to track the grain envelopes. Presently, diffuse-interface-capturing methods [4] are used for the tracking and fixed meshes are employed for the solution of the transport equations. However, diffuse-interface methods can lead to numerical artifacts that are possibly detrimental in describing highly nonlinear phenomena, such as the formation of new branches of the dendrite envelope or the initial stages of growth after nucleation. The use of sharp-interface tracking combined with an adaptive spatial discretization would be advantageous for two reasons: (1) a sharp-interface description would enable us to eliminate the numerical artifacts related to a diffuse-interface; (2) adaptive methods can automatically refine the spatial discretization where required by the form of the solution, typically in the regions of high gradients surrounding the grains, thus increasing the accuracy of the numerical solution.

Because of their flexibility regarding spatial discretization, meshless methods are particularly suitable for the solution of the moving boundary problems encountered in the GEM. The core idea behind the meshless methods is that they operate on scattered nodes [5], in contrast to traditional mesh-based methods that require internodal connectivity. Generally speaking, node placement is simpler than mesh generation, however, not trivial [6]. Dedicated node positioning algorithms for meshless discretization support variable nodal density [7] that can also be applied to parametric surfaces, e.g., a dendrite envelope. In the context of surface reconstruction from a cloud of points, important advances have been just recently made with a two dimensional algorithm using splines [8] and radial basis functions based three-dimensional algorithms [9,10].

In this paper we develop a sharp-interface formulation of the GEM that uses an adaptive meshless discretization method to solve the diffusion in the liquid around the grains. We use the ability of the meshless method to operate on scattered nodes to accurately describe the interface, i.e., the envelope. We demonstrate the ability of the new approach on a two-dimensional $h$-adaptive solution of diffusive growth of dendrites and compare it to a conventional diffuse-interface fixed-grid GEM simulation.

## 2. Mesoscopic Grain Envelope Model for equiaxed isothermal solidification

The Mesoscopic Grain Envelope Model (GEM) [1,2] represents a dendritic grain by its envelope. The growth speed of the envelope, $v_n$, is obtained from the speed of the tips of the dendrite branches (which are not represented in the GEM), $v$, by the relation $v_n = v\cos\theta$, where $\theta$ is the angle between the outward normal to the envelope, $\mathbf{n}$, and the branch growth direction. In the GEM we assume that the dendrite branches grow in predefined growth directions and $\theta$ is the smallest of the angles between $\mathbf{n}$ and the individual directions. In the present study the branches are given $\langle 10 \rangle$ directions, i.e., four possible directions. Note that all equations are written in dimensionless quantities, where the characteristic scales for normalization are $v_{iv}$ for velocity, $D/v_{iv}$ for length and $D/v_{iv}^2$ for time[1]. The tip speed is calculated from a stagnant-film formulation of the 2D *Ivantsov solution* that relates the *Péclet number* of the tip, $Pe$, to $u_\delta$, the normalized concentration difference between the liquid at the tip and that at a finite distance $\delta$ from the tip [11] as

$$u_\delta = \sqrt{\pi Pe}\exp(Pe)\left[\operatorname{erfc}\left(\sqrt{Pe}\right) - \operatorname{erfc}\left(\sqrt{Pe(1+\delta Pe)}\right)\right] \qquad (1)$$

and from a tip selection criterion that reads

$$v = (Pe/Pe_{iv})^2. \qquad (2)$$

Equations (1) and (2) are used to solve for $v$ at any given point on the envelope. The concentration $u_\delta$ is obtained from the concentration field in the liquid around the envelope, which is resolved numerically from the diffusion equation

$$\frac{\partial u}{\partial t} = \nabla^2 u, \qquad (3)$$

---

[1] $v_{iv} = 4\sigma^* D Pe_{iv}^2 / d_0$ is the steady-state velocity of a free (Ivantsov) tip growing into an infinite liquid at supersaturation $\Omega_0$, $\sigma^*$ is the tip selection parameter, $D$ is the diffusion coefficient, $d_0$ is the capillary length, and $Pe_{iv}$ is the Péclet number of the free tip, given by the solution of

$$\Omega_0 = \sqrt{\pi Pe}\exp(Pe_{iv})\operatorname{erfc}\left(\sqrt{Pe_{iv}}\right).$$

subject to boundary conditions $u = 0$ on the envelope and $\mathbf{n} \cdot \nabla u = 0$ on the outer boundary, where $\mathbf{n}$ is the normal vector to the boundary. The initial condition is $u = \Omega_0$, i.e., the liquid is homogeneous and has a dimensionless supersaturation of $\Omega_0$. In this paper we used a stagnant film of $\delta = 1.0$.

### 3. Diffuse-interface GEM on a fixed mesh

The established formulation of the GEM [1,2] uses a phase-field-like interface capturing method [4] to track the grain envelope on a fixed mesh. An indicator function, $\alpha$, is used to distinguish between the interior of the grain ($\alpha > 0.95$) and the liquid ($\alpha < 0.95$). The interface is diffuse since the transition of $\alpha$ between 1 (grain) and 0 (liquid) is smooth with a width $w_\alpha$. The evolution of $\alpha$ is given by

$$\frac{\partial \alpha}{\partial t} + v_n \mathbf{n} \cdot \nabla \alpha = -b \left[ \nabla^2 \alpha - \frac{\alpha(1-\alpha)(1-2\alpha)}{w_\alpha^2} - |\nabla \alpha| \nabla \cdot \left( \frac{\nabla \alpha}{|\nabla \alpha|} \right) \right], \qquad (4)$$

where $v_n$ is the envelope growth speed. The parameters $w_\alpha$ and $b$, as well as the grid size and timestep are selected following the guidelines of Souhar et al. [2]. The advantage of this method is that we can avoid explicit tracking of the envelope and instead simply solve the PDE (4) on a fixed mesh. The drawback is that its accuracy is limited when it needs to describe radii of curvature smaller than around 5 times the mesh spacing [2,4].

### 4. Sharp-interface GEM solved with an adaptive meshless method

In the scope of meshless methods, a linear differential operator $\mathcal{L}$ (for example, $\nabla^2$ from diffusion equation (3)) is approximated from a set of nearby nodes, also referred to as *support nodes* or *stencil*. In this paper, the closest $n$ nodes are chosen as stencil of a node $\mathbf{x}_c$. The approximation of operator $\mathcal{L}$ in node $x_c$ is then sought using an ansatz

$$(\mathcal{L}u)(\mathbf{x}_c) \approx \sum_{i=1}^{n} w_i u_i(\mathbf{x}_i) \qquad (5)$$

for support nodes $\mathbf{x}_i$, support weights $w_i$, an arbitrary function $u$ and nodal values $u_i$. To determine the weights $w_i$, the equality of approximation (5) is enforced for a given set of $m$ basis functions $\{\phi_j\}_{j=1}^{m}$. This allows us to write a linear system

$$\begin{bmatrix} \phi_1(\mathbf{x}_1) & \cdots & \phi_1(\mathbf{x}_n) \\ \vdots & \ddots & \vdots \\ \phi_m(\mathbf{x}_1) & \cdots & \phi_m(\mathbf{x}_n) \end{bmatrix} \begin{bmatrix} w_1 \\ \vdots \\ w_n \end{bmatrix} = \begin{bmatrix} (\mathcal{L}\phi_1)(\mathbf{x}_c) \\ \vdots \\ (\mathcal{L}\phi_m)(\mathbf{x}_c) \end{bmatrix}, \qquad (6)$$

where only the weights $w_i$ are unknown. In this paper monomials with order up to including 2 are used. The support size has been set to $n = 12$, resulting in an overdetermined linear system (6), which is solved using the Householder QR decomposition.

The approximation (5) is general in the sense that it can also be used to interpolate a global solution. In this paper, partition-of-unity (PU) based interpolation has been employed to obtain a global solution approximation [12,13], which is essential for interpolation of concentration or concentration difference (1) at positions that do not correspond to discrete nodal positions.

The stand-alone C++ implementation of the meshless methods, including PU approximators and node positioning algorithms, is available as part of *the Medusa library* [14].

*4.1. Algorithm for the sharp-interface GEM*

The sharp-interface GEM simulation procedure consists of multiple steps shown in figure *1*.

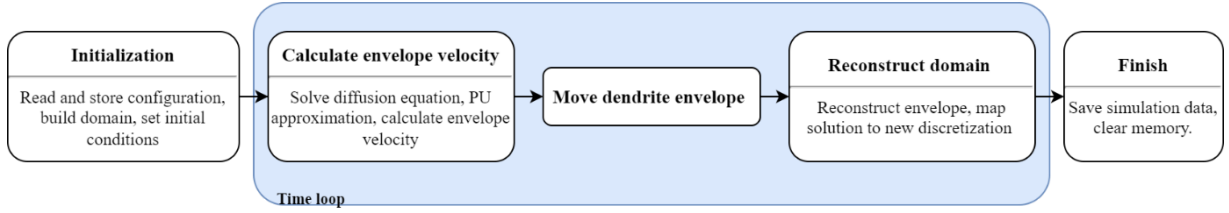

**Figure 1**: Simulation workflow scheme for the meshless sharp-interface GEM.

The initial liquid domain is a difference between the square domain with dimensionless side length $a_m = 20$ and the initial grain envelope represented as a circle with dimensionless radius $r_d = 0.22$. This is the smallest radius of curvature that can be accurately represented with the diffuse-interface method using a grid spacing of 0.05, which is used for comparison in Section 5. The domain is discretized with a linear nodal density distribution between the grain envelope and the outer domain boundary, ranging from $h_d = 0.05$ on the envelope to $h_m = 3h_d$ on the outer boundary. The initial concentration in the liquid is homogeneous at $u = \Omega_0 = 0.18$.

The first step within the explicit time-marching (blue box in **figure 1**) is solving the diffusion equation (3). Afterwards, the envelope velocity is computed for each node on the envelope. In addition to the model parameters, the GEM also requires an accurate concentration field value $u_\delta(\mathbf{x}_i^\delta)$, where $\mathbf{x}_i^\delta = \mathbf{x}_i + \delta \mathbf{n}_i$ for dendrite boundary node $\mathbf{x}_i$, normal $\mathbf{n}_i$ and stagnant-film thickness $\delta$. PU approximation is employed to interpolate a global concentration field $u(\mathbf{x})$ and thus enables us to obtain $u_\delta(\mathbf{x}_i^\delta)$ anywhere in the domain. As a result, the GEM model returns the nodal velocity $v_{tip}$ for the dendrite boundary node $\mathbf{x}_i$ and the boundary node is moved to a new location

$$\mathbf{x}_i' = \mathbf{x}_i + \cos\theta_i v_{tip} dt \mathbf{n}_i \qquad (7)$$

where $\theta_i$ is the minimal angle between the normal $\mathbf{n}_i$ and 4 unit directional vectors, i.e., $\mathbf{e}_x$, $-\mathbf{e}_x$, $\mathbf{e}_y$ and $-\mathbf{e}_y$.

However, moving the envelope boundary nodes in the domain space introduces several problems. Firstly, the internodal distance on the dendritic envelope is generally no longer locally constant and approximately $h$, and secondly, some of the nodes that were originally in the interior of the domain, can fall outside of the computational domain, i.e., in the interior of the dendritic envelope. For the sake of simplicity, we omit these problems by choosing to build a completely new domain discretization with the new dendritic envelope, whose shape is now exactly prescribed with the new boundary nodes from equation (7). To accurately build the domain discretization and assure minimal distortion to the dendritic envelope, a surface reconstruction algorithm is required. In this paper we reconstruct the two-dimensional envelope shape using splines, as in [8]. After the discretization is completed, the concentration values in the new nodal positions are obtained by employing the PU interpolation constructed over the old domain discretization.

The entire process, i.e., the blue rectangle on figure 1, is then repeated, until the total simulation time $t_{tot} = 12$ is reached.

## 5. Results

The initial domain has been filled with 21945 nodes and computational time of approximately 35 hours on a computer with `AMD EPYC 7702 64-Core` processor with `512GB of DDR4` memory for simulation time $t_{tot} = 12$ with timestep $dt = 10^{-4}$. The code[2] was compiled using `g++ (GCC) 9.3.0 for Linux` with `-O3 -std=c++17 -fopenmp -DNDEBUG` flags.

A timelapse of the growing grain is shown in figure 2.

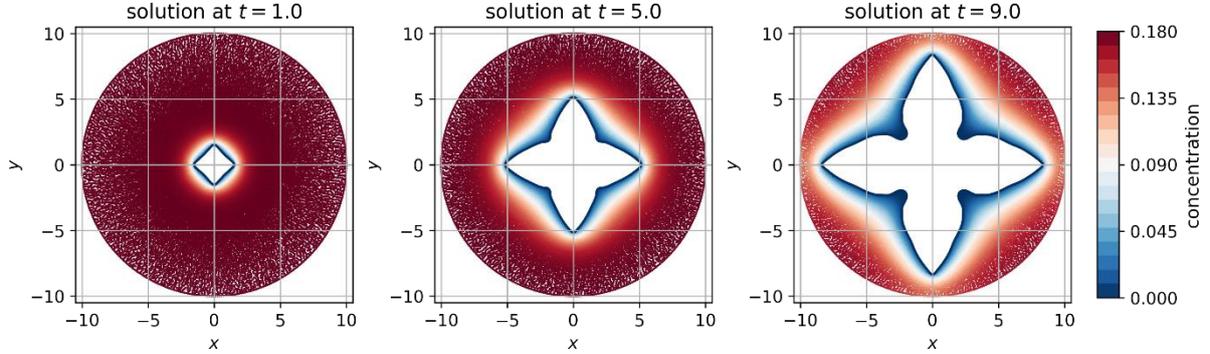

**Figure 2:** Dendrite growth timelapse.

As seen in figure *2*, we are able to obtain a reasonable solution using the sharp-interface meshless method. To assess its quality, we compare it to the diffuse-interface fixed-mesh solution obtained with identical physical and model parameters and a mesh spacing of 0.05. The primary tip speed as a function of the primary tip position is compared in figure *3* (left) and the grain envelope shape at $t = 6$ in figure *3* (right).

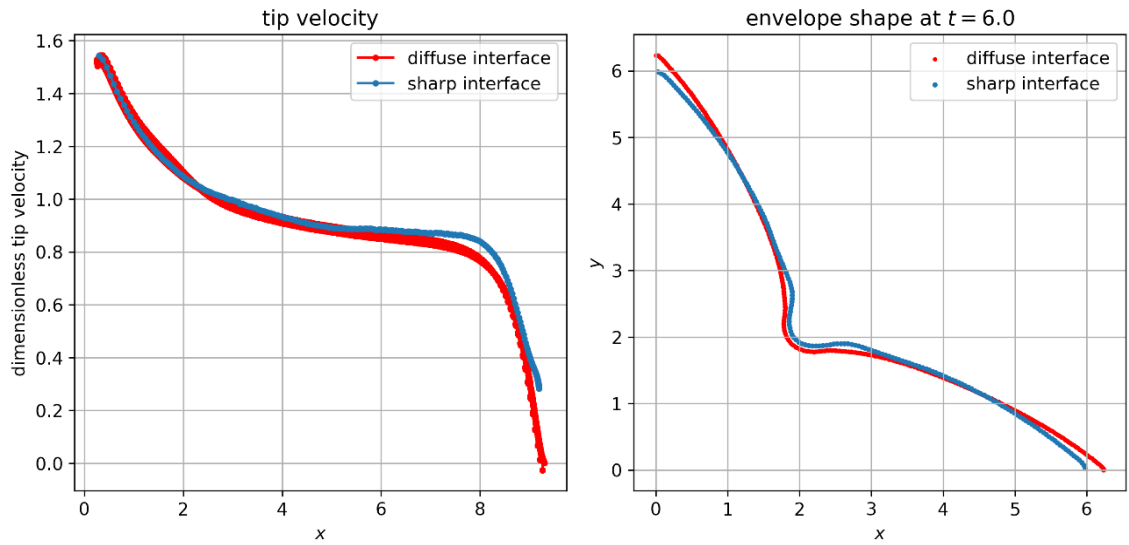

**Figure 3:** Dimensionless tip velocity vs. tip position (left) and dendrite envelope shape at selected simulation time (right).

The evolution of the tip velocity shows the characteristic initial transient form a high initial velocity towards a velocity of 1. In the present case the tip velocity does not reach a steady state since the diffusion field interacts with the outer domain boundary before steady state is reached. The final transient consists of a slowdown of the tip that finally stops at a distance $\delta$ from the boundary, i.e., $x_{tip} \approx 9$.

---

[2] Source code available at: https://gitlab.com/e62Lab/public/2022_cp_icasp_dendrite_growth under tag *v1.0*.

## 6. Conclusions

We have shown that the *h*-adaptive meshless method can be successfully employed to simulate the growth of a dendrite envelope during solidification, using a sharp-interface representation of the envelope. The meshless numerical solution has been compared to the conventional diffuse-interface fixed-grid GEM solution in terms of maximum tip velocity for all the simulation steps and in terms of tip shape at a selected time step. We have shown that both solutions obtained by different numerical methods are comparable. Future work could include extension to the more realistic three-dimensional dendrites and full exploitation of parallelism to reduce the computational times of computationally more expensive meshless methods.

**Acknowledgments.** The authors would like to acknowledge the financial support of the Slovenian Research Agency (ARRS) research core funding No. P2-0095 and the World Federation of Scientists.

**References**
[1] Steinbach I, Diepers H-J and Beckermann C 2005 Transient growth and interaction of equiaxed dendrites *J. Cryst. Growth* **275** 624–38
[2] Souhar Y, De Felice V F, Beckermann C, Combeau H and Založnik M 2016 Three-dimensional mesoscopic modeling of equiaxed dendritic solidification of a binary alloy *Comput. Mater. Sci.* **112** 304–17
[3] Tourret D, Sturz L, Viardin A and Založnik M 2020 Comparing mesoscopic models for dendritic growth *IOP Conference Series: Materials Science and Engineering* vol 861 (IOP Publishing) p 012002
[4] Sun Y and Beckermann C 2007 Sharp interface tracking using the phase-field equation *J. Comput. Phys.* **220** 626–53
[5] Wang H and Qin Q-H 2019 Overview of meshless methods *Methods of Fundamental Solutions in Solid Mechanics* (Elsevier) pp 3–51
[6] Shankar V, Wright G B, Kirby R M and Fogelson A L 2015 A Radial Basis Function (RBF)-Finite Difference (FD) Method for Diffusion and Reaction–Diffusion Equations on Surfaces *J. Sci. Comput.* **63** 745–68
[7] Slak J and Kosec G 2019 On Generation of Node Distributions for Meshless PDE Discretizations *SIAM J. Sci. Comput.* **41** A3202–29
[8] Jancic M, Cvrtila V and Kosec G 2021 Discretized Boundary Surface Reconstruction *2021 44th International Convention on Information, Communication and Electronic Technology (MIPRO)* 2021 44th International Convention on Information, Communication and Electronic Technology (MIPRO) (Opatija, Croatia: IEEE) pp 278–83
[9] Drake K P, Fuselier E J and Wright G B 2022 Implicit Surface Reconstruction with a Curl-free Radial Basis Function Partition of Unity Method *ArXiv210105940 Cs Math*
[10] Liu S, Liu T, Hu L, Shang Y and Liu X 2021 Variational progressive-iterative approximation for RBF-based surface reconstruction *Vis. Comput.* **37** 2485–97
[11] Cantor B and Vogel A 1977 Dendritic solidification and fluid flow *J. Cryst. Growth* **41** 109–23
[12] Slak J 2021 Partition-of-Unity Based Error Indicator for Local Collocation Meshless Methods *2021 44th International Convention on Information, Communication and Electronic Technology (MIPRO)* 2021 44th International Convention on Information, Communication and Electronic Technology (MIPRO) (Opatija, Croatia: IEEE) pp 254–8
[13] De Marchi S, Martínez A and Perracchione E 2019 Fast and stable rational RBF-based partition of unity interpolation *J. Comput. Appl. Math.* **349** 331–43
[14] Slak J and Kosec G 2021 Medusa: A C++ Library for Solving PDEs Using Strong Form Mesh-free Methods *ACM Trans. Math. Softw.* **47** 1–25